\documentclass{article}
\usepackage{graphicx} % Required for inserting images
\usepackage{hyperref}
\usepackage[ruled, vlined, linesnumbered]{algorithm2e}
\usepackage{amsmath}
\usepackage{multirow} % For multi-row and multi-column cells
\usepackage{array}    % For better table control
\usepackage{changepage}
\usepackage{csquotes}

\title{\textbf{Wave Function Collapse Set Covering and the Hill Climbing Algorithm: A New, Fast Heuristic \& Metaheuristic Pairing for the Minimum Set Cover Problem}}
\author{David Oprea, David Perkins}
\date{August 2024}

\begin{document}

\maketitle

\begin{abstract}
In this paper, we present a new heuristic that focuses on the optimization problem for the Minimum Set Cover Problem. Our new heuristic involves using Wave Function Collapse and is called Wave Function Collapse Set Covering (WFC-SC). This algorithm goes through observation, propagation, and collapsing, which will be explained more later on in this paper. We optimize this algorithm to quickly find optimal coverings that are close to the minimum needed. To further optimize this algorithm, we pair it with the Hill Climbing metaheuristic to help WFC-SC find a better solution than its “local optimum." We benchmark our algorithm using the well known OR library from Brunel University against other proficient algorithms. We find that our algorithm has a better balance between both optimality and time.
\end{abstract}

\section{\textbf{Introduction}}

For a given set system on a universe of items and a collection of a set of items, the Minimum Set Cover Problem (MSCP) \cite{JOHNSON1974256} for optimization finds the minimum number of sets that covers the whole universe. This problem was proven to be NP Hard by Karp \cite{Karp1972}. Exact algorithms for this problem exist like backtracking, integer linear programming (ILP), and branch and cut \cite{caprara2000algorithms}, but once the sets get large enough, exact algorithms become too slow and resource-expensive to use. This is why researchers tend to focus on heuristics and metaheuristics for solving problems like these. Examples include: greedy algorithms, genetic algorithms, local search, and ant colony optimization \cite{metaheuristics}. There has also been recent research to further improve upon these heuristics. For example, the Big Greedy heuristic \cite{Chandu2015ImprovedGA} attempts to improve upon the greedy heuristic and incorporates local search and Tabu Search \cite{TabuSearch} to reduce the number of sets in a solution found by the greedy heuristic initially. 

The problem has numerous applications in different areas of study and industrial applications. The applications include multiple sequence alignments for computational biochemistry, manufacturing, network security, service planning and location problems. 

These heuristics are able to find a solution quickly due to their simplicity, but one major problem with them is that their solutions could be more optimal. One thing most of them miss is a metaheuristic to help further guide them to a more optimal solution instead of sticking to a local optimum. Our first goal of our research paper is to present a new heuristic, WFC-SC that is not only just as fast as these other effective heuristics, but it can also be just as optimal. Our second goal is to show that our heuristic, paired with the Hill Climbing metaheuristic will be able to perform better with little time loss. We hope that our research could further inspire the use of heuristic and metaheuristic pairs toward solving other NP problems.  

The structure of the paper is as follows: in Section 2, we discuss the formal definition of the Minimum Set Cover Problem; in Section 3, we outline previous fast heuristics for the MSCP that we will compare our heuristic to; in Section 4, we present the formal algorithm for WFC-SC; in Section 5, we showcase WFC-SC with Hill Climbing and how to perform it; in Section 6, we generate experimental results and compare such results to existing literature; in Section 7, we conclude with a discussion on the uses of WFC-SC and possible future research regarding WFC-SC.

\section{\textbf{The Minimum Set Cover Problem}}

Given a set of elements {1, 2, …, n} (called the universe) and a collection S of m subsets whose union equals the universe, the set cover problem is to identify the smallest sub-collection of S whose union equals the universe.

More formally, given a universe $U$ and a family $S$ of subsets of $U$, a set cover is a subfamily $C$ of $S$ whose sets' union is $U$.

\begin{center}
\includegraphics[scale = 0.5]{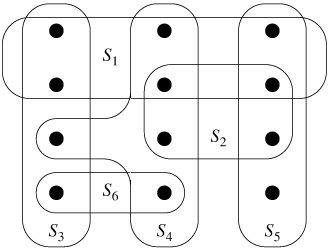}
\end{center}

In the image above, we have 12 dots which represent the elements that make up the universe $U$, and we have multiple subsets where each element is being covered by at least one of them. 

In this case, the most optimal answer would be to take $S_3$, $S_4$, and $S_5$ (all other valid answer have more than three sets). The goal of the minimum set cover problem is to try to find this optimal number for any possible universe and subsets. This will be the main focus of our heuristic and metaheuristic pairing throughout the rest of the paper. 

\section{\textbf{Heuristic Approaches}}
\label{Heuristics}

In this section, we go through a couple of other heuristic approaches that are popular and effective at finding approximate answers to this problem. These algorithms are \textit{Greedy}, \textit{Big Greedy}, and \textit{Tabu Search}. Section 6 compares the performance of WFC-SC + Hill Climbing to the performance of the previously mentioned algorithms. 

\subsection{Greedy}
\label{GreedySec}

The greedy algorithm \cite{Greedy} is the simplest heuristic for the set cover problem. All it does is insert a set that adds the most new elements to its answer until its answer covers all elements in the universe $U$. It answers the question in approximately O($\log N$) time and usually gives a solution close to the minimum. 

\begin{algorithm}[H]
    \SetAlgoLined
    \SetKwInOut{KwIn}{Input}
    \SetKwInOut{KwOut}{Output}
    \KwIn{The universe $U$, a family $S$ of subsets that covers $U$, and an empty set $I$}
    \KwOut{A set cover for $U$}
    \While{$\bigcup_{A \in I}{A} \neq U$}{
        $s \leftarrow$ a set in $S$ with the most elements not yet in $\bigcup_{A \in I}{A}$\;
        $push(I, s)$ \;
        $remove(S, s)$ \;
    }
    return $I$ \;
    \caption{Greedy}
    \label{Greedy}
\end{algorithm}

\subsection{Big Greedy}

Big Greedy \cite{Chandu2015ImprovedGA} is a more advanced version of the greedy algorithm presented previously. The only difference this heuristic has is a step-size parameter $p$ which controls the number of sets that will be added to the answer next. 

\begin{algorithm}[H]
    \SetAlgoLined
    \SetKwInOut{KwIn}{Input}
    \SetKwInOut{KwOut}{Output}
    \KwIn{The universe $U$, a family $S$ of subsets that covers $U$, the step-size parameter $p$, and an empty set $I$}
    \KwOut{A set cover for $U$}
    \While{$\bigcup_{A \in I}{A} \neq U$}{
        \eIf{$p > |S|$}{
            $P \leftarrow S$\;
        }{
            $P \leftarrow$ the $p$ sets currently in $S$ whose union has the most elements not yet in $\bigcup_{A \in I}{A}$\;
        }
        \For{each set $s$ in $P$}{
            $push(I, s)$ \;
            $remove(S, s)$ \;
        }
    }
    return $I$ \;
    \caption{Big Greedy}
    \label{BigGreedy}
\end{algorithm}

\subsection{Tabu Search}
\label{TabuSec}

Tabu Search \cite{TabuSearch} uses a local search and a tabu list to try to find the minimum set covering. The local search involves having a set $ans$ as the current possible answer and neighbors being the sets in $S$. One can either add this set if it is not in ans or remove it otherwise. The tabu list, with a length of $t$, is used to prevent repetitive steps from occurring to help the search find a more optimal answer. Also, the researchers include an upper bound parameter $u$ which is the max number of sets that ans can have. This is used to make sure their search is always trying to find a better answer. 

To score the fitness of ans, the researchers use a very simple fitness function:

\vspace{1em}

\begin{center}
\textbf{$Fitness = NrOfUncoveredElements + |ans|$}.
\end{center}

\vspace{1em}

\noindent Where $NrOfUncoveredElements$ is the number of elements in the universe $U$ that ans still does not contain and $|ans|$ is the number of sets currently in ans.

\begin{algorithm}[H]
    \SetAlgoLined
    \SetKwInOut{KwIn}{Input}
    \SetKwInOut{KwOut}{Output}
    \KwIn{The universe $U$,  a family $S$ of subsets that covers $U$, the tabu list $T$, the length of the tabu list $t$, the number of iterations $n$, and a set $I$ which is the answer from the greedy algorithm (see Section \ref{GreedySec})}
    \KwOut{A set cover for $U$}
    $best \leftarrow$ $I$ \;
    $u \leftarrow  |I| - 1$\;
    \For{i in 0 to $n$}{
        $s \leftarrow$ the next set in $S$ to be added or removed from $I$, where it is not in $T$ and the new length is not greater than $u$, with the best next $fitness$\;
        \eIf{$s$ in $I$}{
            $remove(I, s)$ \;
        }{
            $push(I, s)$ \;
        }
        $push(T, s)$ \;
        \If{$|T| > t$}{
            remove first element from $T$ \;
        }
        \If{$\bigcup_{A \in I}{A} = U$}{
            $best \leftarrow$ $I$ \;
            $u \leftarrow |best| - 1$\;
        }
    }
    return $best$ \;
    \caption{Tabu Search}
    \label{TabuSearch}
\end{algorithm}

\section{Wave Function Collapse}

Before showing the pseudocode, we will give a simple example of how WFC-SC actually works. 

Say we are given a universe
\[
U = \{1, 2, 3, 4, 5\}
\]
and we have some number of subsets whose union equals $U$, like
\[
S_1 = \{1, 4\}, \quad
S_2 = \{1, 2, 3\}, \quad
S_3 = \{2, 3, 5\}, \quad
S_4 = \{3, 4\}. \quad
\]

We define an \textit{active element} to be an element that we have not yet included in our answer and we define an \textit{active set} to be a set we have not yet removed from $S$ or included in our answer. Initially, every set and element is active. We define an active set's number of \textit{conflicts} to be the sum of the number of times each of its active elements appears in all active sets other than itself. We define an active set's \textit{entropy} to be the number of active elements in that set. So, the conflicts and entropy for each of the sets in the example is:

\begin{align*}
S_1 &= \{1, \; 4\}; conflicts = 2, entropy = 2,\\
S_2 &= \{1, \; 2, \; 3\}; conflicts = 4, entropy = 3, \\
S_3 &= \{2, \; 3, \; 5\}; conflicts = 3, entropy = 3, \\
S_4 &= \{3, \; 4\}; conflicts = 3, entropy = 2.
\end{align*}

First, we check for each active element if it is only contained in one set. In this case, 5 is only contained in $S_3$, so $S_3$ has to be added to the answer. Now since we added it to the answer, the conflicts and entropy for the other sets now change. 

\begin{align*}
S_1 &= \{1, \; 4\}; conflicts = 2, entropy = 2,\\
S_2 &= \{1, \; 2, \; 3\}; conflicts = 1, entropy = 1, \\
S_4 &= \{3, \; 4\}; conflicts = 1, entropy = 1.
\end{align*}

Remember that conflicts and entropy are only counted for active elements. So, we would not consider the conflicts or entropy added by elements 2 and 3 since they are part of $S_3$ and were added to the answer. 

Now we remove the set with the worst \textit{average conflicts} which is equal to a set's conflicts divided by its entropy, breaking ties using the lowest entropy and then choosing randomly. In this case, we can remove any of the sets since all their average conflicts are one. Let us remove $S_4$. 

\begin{align*}
S_1 &= \{1, \; 4\}; conflicts = 1, entropy = 2,\\
S_2 &= \{1, \; 2, \; 3\}; conflicts = 1, entropy = 1. \\
\end{align*}

We then check the active elements that were in $S_4$ and see if they are now contained by only one set. In this case 4 is now only contained by $S_1$, so we add it to the answer. Now we have: 

\begin{align*}
S_2 &= \{1, \; 2, \; 3\}; conflicts = 0, entropy = 0. \\
\end{align*}

We remove any sets that contain an entropy of 0 because they cannot add any new elements to the answer, so now we finished checking all of the subsets. Our answer is

\begin{align*}
answer &= \{S_1, \; S_3\}.
\end{align*}

One last definition is needed before we can understand the pseudocode presented later. We define $contains(x)$ as a function that returns a list of the active sets in the set of subsets $S$ that contain the element $x$. So if we were to go back to our previous example and look at the subsets:

\[
S_1 = \{1, 4\}, \quad
S_2 = \{1, 2, 3\}, \quad
S_3 = \{2, 3, 5\}, \quad
S_4 = \{3, 4\}, \quad
\]

\noindent and asked what $contains(3)$ would be, it would return: 

\begin{align*}
\{S_2, \; S_3\, \; S_4\},
\end{align*}

\noindent since each of those contain 3 and are initially active. Though if one of those sets were to become inactive, then the value returned would be different. 

Our pseudocode for our heuristic, WFC-SC, is shown below and our metaheuristic, Hill Climbing, is shown later on in Section \ref{MetaHeuristic}. 

\begin{algorithm}[H]
    \SetAlgoLined
    \SetKwInOut{KwIn}{Input}
    \SetKwInOut{KwOut}{Output}
    \KwIn{The universe $U$, a family $S$ of subsets that covers $U$}
    \KwOut{A set cover for $U$}
    \For{each value $x$ in $U$}{
        \If{only one active set $a$ contains $x$}{
            $collapse(S, a)$
        }
    }
    \While{$S$ is not empty}{
        $s \leftarrow$ $observe(S)$ \;
        $propagate(S, s)$ \;
    }
    \caption{Set Covering via WFC-SC}
    \label{WFC-SC}
\end{algorithm}  

Our observe function scores sets based on their average conflicts which we defined previously in our example. 

\begin{algorithm}[H]
    \SetAlgoLined
    \SetKwInOut{KwIn}{Input}
    \SetKwInOut{KwOut}{Output}
    \KwIn{A set $S$ with the current \textit{active sets}}
    \KwOut{The set with the highest \textit{average conflicts}}
    return the set in $S$ with the highest \textit{average conflicts}, breaking ties with lowest $entropy$  \;
    \caption{Observe}
    \label{Observe}
\end{algorithm}

The propagate function removes the set returned by observe from the set $S$ of possible sets. It then goes through the active elements in that set and reduces the conflicts from the other sets by one for each of those active elements it contains. Finally, it checks each active element to see if only one active set contains it. If this is the case, that active set is put into our collapse function. 

\begin{algorithm}[H]
    \SetAlgoLined
    \SetKwInOut{KwIn}{Input}
    \KwIn{A set $S$ with the current \textit{active sets} and the set $s$ to be removed from $S$}
    $remove(S, s)$ \;
    \For{each \textit{active number} $x$ in set $s$}{
        \For{each set $t$ in $contains(x)$}{
            $t.conflicts \leftarrow t.conflicts - 1$ \;
        }
        \If{only one active set $a$ contains $x$}{
            $collapse(S, a)$ \;
        }
    }
    \caption{Propagate}
    \label{Propagate}
\end{algorithm}

We define one more element, $ans$, as a list containing the optimal amount of subsets needed, whose union creates the universe $U$. The collapse function adds the set $s$ into ans. It then goes through the set's active elements, removes them for consideration, and updates the necessary active sets. Finally, it removes the set $s$ and other active sets whose entropy has become zero from $S$. 

\begin{algorithm}[H]
    \SetAlgoLined
    \SetKwInOut{KwIn}{Input}
    \KwIn{A set $S$ with the current \textit{active sets} and the set $s$ to be added to $ans$}
    $push(ans, s)$ \;
    \For{each \textit{active number} $x$ in set $s$}{
        \For{each set $t$ in contains($x$)}{
            $t.conflicts \leftarrow t.conflicts - |contains(x)|$ \;
            $t.entropy \leftarrow t.entropy - 1$ \;
            \If{$t.entropy = 0$}{
                $remove(S, t)$ \;
            }
        }
    }
    $remove(S, s)$ \;
    \caption{Collapse}
    \label{Collapse}
\end{algorithm}

\vspace{1em}
\hrule 
\vspace{1em}

\section{Metaheuristic (Hill Climbing)}
\label{MetaHeuristic}

Sometimes the Wave Function Collapse algorithm will reach a local optimum, but it is not able to find the global optimum. This is where our metaheuristic, the Hill Climbing Algorithm \cite{hillClimbing}, comes in. 

Firstly, we make a slight adjustment to our observe function, which changes how subsets are scored for our WFC-SC function that we will call $score$. The value of score is equal to

\vspace{1em}

\begin{center}
\textbf{$s.conflicts^{c} / s.entropy^{exp}$}
\end{center}

\vspace{1em}

\noindent where $exp$ will be the value the Hill Climbing metaheuristic will adjust and $c$ is a constant value. For our testing, we made c equal to 0.9 for all files except the scp60 files, where we made c equal to 0.7 for better performance of our algorithm. 

\begin{algorithm}[H]
    \SetAlgoLined
    \SetKwInOut{KwIn}{Input}
    \SetKwInOut{KwOut}{Output}
    \KwIn{A set $S$ with the current \textit{active sets} and the exponent factor $exp$}
    \KwOut{The set with the highest $score$}
    return the set in $S$ with the highest $score$, breaking ties with lowest $entropy$  \;
    \caption{Observe (For Hill Climbing)}
    \label{ObserveHill}
\end{algorithm}

We also need to slightly redefine our WFC-SC function reference in Algorithm \ref{WFC-SC} to fit with our metaheuristic. We change our algorithm's inputs to the universe $U$, a family $S$ of subsets that covers $U$, and a new $exp$ parameter which will be used in the new observe function. We also change line 7 from $s \leftarrow observe(S)$ to $s \leftarrow observe(S, exp)$. These are all the changes we need to implement our metaheuristic. 

Again, the Hill Climbing metaheuristic's goal is to adjust the value $exp$ so that our wave function collapse can find the most optimal set covering. With every iteration, the algorithm tests the current value of exp using WFC-SC and looks at the length of the new answer. If that answer is better than the previous one, then the algorithm will keep moving exp in the same direction (increase or decrease). Otherwise, it will change the direction. After each iteration, the value by which exp changes decreases slightly to help it land into a suitable answer. This goes on for a number of steps, which we assign the variable to the variable $n$.  

\begin{algorithm}[H]
    \SetAlgoLined
    \SetKwInOut{KwIn}{Input}
    \SetKwInOut{KwOut}{Output}
    \KwIn{The universe $U$, a family $S$ of subsets that covers $U$, the number of iterations $n$, the initial value $exp_i$ of $exp$, and the temperature $T$ for adjusting $exp$}
    \KwOut{A set cover for $U$}
    $exp \leftarrow exp_i$ \;
    $best \leftarrow$ The solution from WFC-SC when $exp = 1$ \;
    $prev \leftarrow best$ \;
    $add \leftarrow 0$ \; 
    \For{i in 0 to $n$}{
        $s \leftarrow$ WFC-SC($U$, $S$, $exp$) \;
        \eIf{$|prev| > |s|$}{
            \eIf{$add > 0$}{
                $add \leftarrow T * exp$ \;
            }{
                $add \leftarrow -T * exp$ \;
            }
            \If{$|best| > |s|$}{
                $best \leftarrow s$ \;
            }
        }{
            \eIf{$add > 0$}{
                $add \leftarrow -T * exp$ \;
            }{
                $add \leftarrow T * exp$ \;
            }
        }
        $prev \leftarrow s$ \;
        $T \leftarrow T * 0.99$ \;
        $exp \leftarrow exp + add$ \;
    }
    return $best$ \;
    \caption{Hill Climbing}
    \label{Hill Climbing}
\end{algorithm}

\section{Experimental Results}

Here, we compare WFC-SC’s performance to the algorithms described earlier (in Section \ref{Heuristics}): Greedy, Big Greedy, and Tabu Search. 

\subsection{Experimental Environment}

We implemented this algorithm in Julia 1.10, on a 3.0 GHz 11th Gen Intel Core i7 processor. To test the performance of WFC-SC, we used a popular set cover dataset, the OR benchmark from Brunel University which is introduced in \cite{ORLibrary}. It contains both testing sets for the uni-cost set cover problem (which we worked on), and the multi-cost set cover problem. The OR benchmark has been used for many of the algorithms which have been presented in this paper. It contains both the testing data and answers for the best answers for each of its test files. 

The test files we use for testing were presented in the following papers: the scp40, scp50, and scp60 files were originally from the Balas and Ho set covering paper \cite{Balas1980}, the scpd and scpe files are the test problem sets from J.E Beasley's paper on set covering \cite{BEASLEY198785}, and the scpclr files were contributed to the benchmark by A. Wool.

In Table \ref{Results}, you will see the results for each of the algorithms. For each of them we measure two variables: the number of sets in their answer, $k$, and the mean time they took to get that answer. Also, for Tabu Search we measure an additional variable, $Best$, which is the best answer Tabu Search found during its iterations. Note the differences in the metric prefixes for the mean times of each algorithm which was used to make the formatting more appealing. For Tabu Search, we ran the algorithm 25 times for each file since its tiebreaker method involves randomness (this is shown in Algorithm \ref{TabuSearch}). 

\vspace{1em}

\begin{table}[h!]
    \caption{Comparison Between WFC-SC + Hill Climbing and other Heuristics in the OR Benchmark} % Title of the table
    \centering % Center the table
    \resizebox{\textwidth}{!}{ % Resize the table to fit the text width
    \begin{tabular}{|c|c|c|c|c|c|c|c|c|c|}
        \hline
        \multirow{2}{*}{Test Set} & \multicolumn{2}{c|}{Greedy} & \multicolumn{2}{c|}{Big Greedy \cite{Chandu2015ImprovedGA}} & \multicolumn{3}{c|}{Tabu Search \cite{TabuSearch}} & \multicolumn{2}{c|}{WFC-SC + Hill Climbing} \\ \cline{2-10}
         & $k$ & Mean Time(ns) & $k$ & Mean Time(ns) & $k$ & Best & Mean Time(s) & $k$ & Mean Time(ms) \\ \hline
        scp41.txt & 41 & 149.40 & 42 & 168.20 & 40.36 & 39 & 3.21 & 41 & 25.25 \\ \hline
        scp42.txt & 41 & 151.86 & 41 & 162.82 & 38.32 & 37 & 6.79 & 39 & 26.34 \\ \hline
        scp43.txt & 43 & 159.47 & 43 & 155.00 & 39.92 & 38 & 6.51 & 41 & 24.56 \\ \hline
        scp44.txt & 44 & 159.78 & 43 & 64.15 & 40.96 & 40 & 2.53 & 41 & 26.30 \\ \hline
        scp51.txt & 37 & 67.96 & 38 & 146.02 & 36.56 & 36 & 4.23 & 38 & 12.18 \\ \hline
        scp52.txt & 38 & 83.72 & 38 & 70.25 & 36.4 & 35 & 5.32 & 38 & 22.91 \\ \hline
        scp53.txt & 37 & 161.05 & 39 & 161.03 & 36.28 & 35 & 4.09 & 38 & 23.96 \\ \hline
        scp54.txt & 39 & 83.72 & 38 & 68.50 & 36.16 & 35 & 3.92 & 39 & 12.52 \\ \hline
        scp61.txt & 23 & 63.70 & 24 & 82.09 & 22.04 & 21 & 2.27 & 23 & 9.41 \\ \hline
        scp62.txt & 23 & 71.83 & 24 & 143.62 & 22.44 & 22 & 1.92 & 23 & 10.30 \\ \hline
        scp63.txt & 23 & 92.55 & 23 & 152.69 & 22.20 & 22 & 1.73 & 23 & 11.25 \\ \hline
        scpe1.txt & 5 & 151.33 & 6 & 65.56 & 5.00 & 5 & 2.53 & 5 & 5.41 \\ \hline
        scpe2.txt & 5 & 75.31 & 6 & 70.18 & 5.00 & 5 & 3.90 & 5 & 5.79 \\ \hline
        scpe3.txt & 5 & 64.45 & 6 & 63.57 & 5.00 & 5 & 1.70 & 5 & 2.41 \\ \hline
        scpclr11.txt & 33 & 163.40 & 31 & 192.73 & 27.28 & 26 & 3.12 & 30 & 34.81 \\ \hline
        scpclr12.txt & 30 & 72.68 & 30 & 62.99 & 28.36 & 24 & 8.44 & 30 & 60.60 \\ \hline
        scpclr13.txt & 32 & 86.57 & 31 & 72.77 & 28.36 & 25 & 11.35 & 32 & 128.10 \\ \hline
        % Add more rows as needed
    \end{tabular}
    }
    \label{Results}
\end{table}

\begin{table}[h!]
\centering
\begin{tabular}{|l|c|c|}
\hline
\textbf{Test File} & \textbf{Best} & \textbf{WFC + Hill Climbing + Tabu} \\
\hline
s12.txt (41) & 38 & 40 \\
s13.txt (42) & 37 & 39 \\
s14.txt (43) & 38 & 39 \\
s15.txt (44) & 38 & 42 \\
s16.txt (51) & 35 & 37 \\
s17.txt (52) & 35 & 36 \\
s18.txt (53) & 34 & 38 \\
s19.txt (54) & 34 & 36 \\
s20.txt (e1) & 5  & 5  \\
s21.txt (e2) & 5  & 5  \\
s22.txt (e3) & 5  & 5  \\
s23.txt (61) & 21 & 22 \\
s24.txt (62) & 20 & 23 \\
s25.txt (63) & 21 & 23 \\
\hline
\end{tabular}
\caption{Comparison of best known results with WFC + Hill Climbing + Tabu search}
\end{table}

\vspace{1em}

\subsection{Results}

The data for the scp40 files contained 1000 sets with the elements being integers from 1-200. Based on our experimental data for the scp40 files, we find that Tabu Search outperforms WFC-SC by less than one set on average and that its best is about two less sets than WFC-SC on average. Furthermore, we find that WFC-SC outperforms Greedy and Big Greedy by about two sets on average for these files. 

The data for the scp50 files contained 2000 sets with the elements being integers from 1-200. Based on our experimental data for the scp50 files, we find that Tabu Search outperforms WFC-SC by about two sets on average and that its best is about three sets less than WFC-SC. We also find that WFC-SC performs slightly worse than Greedy and that Big Greedy and WFC-SC are on par with each other. 

The data for the scp60 files contained 1000 sets with the elements being integers from 1-200. When looking at the data throughout the scp60 files we find that WFC-SC is on par with Greedy, Tabu Search performs better than WFC-SC by less than one set on average and its best is about one set less than WFC-SC. Big Greedy performs worse than WFC-SC by one set on two of the files. 

The data for the scpe files contained 500 sets with the elements being integers from 1-50. We find that Greedy, Tabu Search, and WFC-SC all perform the same with Big Greedy performing one set worse for all of these files. 

The data for the scpclr files is different for all of them. The scpclr11 file contained 210 sets with the elements being integers from 1-511. The scpclr12 file contained 330 sets with the elements being integers from 1-1023. The scpclr13 file contained 495 sets with the elements being integers from 1-2047. Throughout the scpclr files, we find that Tabu Search outperforms WFC-SC by about three sets on average and that its best is about five sets less than WFC-SC on average. Furthermore, we find that WFC-SC outperforms Greedy and Big Greedy by about one set on average for these files. 

Throughout all the files, Tabu Search is clearly significantly slower than WFC-SC, whereas the greedy heuristics are consistently a lot faster and on par with one another in terms of speed and optimality, but WFC-SC on average performs better than both of these heuristics. 

\section{Conclusion and Further Research}

In this paper, we presented a new method for solving the Set Covering Problem: Wave Function Collapse (WFC-SC) with Hill Climbing. WFC-SC goes through three main stages: observe, propagate, and collapse in order to find an optimal set covering for some test. We then improved upon this by augmenting the scoring function such that WFC-SC's observe function uses through our hill climbing metaheuristic and an exp parameter.

Even though our algorithm does not consistently outperform the other algorithms we compared to, it does perform decently well in both optimality and time. Though Greedy and Big Greedy are significantly faster, WFC-SC combined with Hill Climbing tends to do better than Greedy and Big Greedy. And although Tabu Search does perform better on average, WFC-SC + Hill Climbing is comparatively a lot faster (about 500 times faster) compared to Tabu Search. All in all, our algorithm balances both speed and optimally using a simple method. 

For future work and research, this format of using Wave Function Collapse and a metaheuristic pairing could work well for other NP hard problems. With our limited time, we could undoubtedly have optimized our observe function's scoring or metaheuristic or have tested on more of the test files presented in the OR library for more accurate results. From this, further work in for our research could involve fine-tuning our WFC-SC functions or testing our a different metaheuristic with it to see if it can be improved. We also did not consider using this algorithm for the multi cost set covering problem, but with minor modifications to our code it could also solve for that problem. 

% This next line includes the bibliography in your paper:

\end{document}